\documentclass{amsart}
\usepackage{amsmath, amssymb}
\newtheorem{Theorem}{Theorem}[section]
\newtheorem{Cor}[Theorem]{Corollary}
\newtheorem{Lemma}[Theorem]{Lemma}

\newtheorem{Proposition}[Theorem]{Proposition}
\theoremstyle{definition}
\newtheorem{Definition}[Theorem]{Definition}
\theoremstyle{remark}

\numberwithin{equation}{section}
\newcommand{\R}{\mathbb R}
\newcommand{\N}{\mathbb N}

\newcommand{\C}{\mathcal{C}}
\newcommand{\F}{\mathcal{F}}
\newcommand{\p}{\mathcal{P}}
\newcommand{\U}{\mathcal{U}}

\newcommand{\rel}{\mathcal{R}}

\newcommand{\B}{\mathcal{B}}

\begin{document}

\title{From configurations to branched configurations and beyond}%
\author{Jean-Pierre Magnot}
\address{ Lyc\'ee Jeanne dArc \\ 30 avenue de Grande Bretagne\\
F-63000 Clermont-Ferrand}
\email{jean-pierr.magnot@ac-clermont.fr}%

%\date{10/11/2012}
\begin{abstract}
We propose here a geometric and topological setting for the study of branching effects arising in various fields of research, e.g. in statistical mechanics and turbulence theory. We describe various aspects  
that appear key points to us, and finish with a limit of such a construction which stand in the dynamics on probability spaces where it seems that branching effects can be fully studied without any adaptation of the framework.      
\end{abstract}

\maketitle

                     \tableofcontents

\section*{Introduction}
Finite and infinite configuration spaces are rather old topics, 
see e.g. \cite{Str, Ism} and the references cited therein,
that had many applications in various settings in mathematical physics and representation theory. 
 
More recently, several papers, including \cite{FKK0,FKK1,FKK2,FKK3} showed how these topics could be applied in various disciplines:ecology, financial markets, and so on. This large spectrum of applications principally comes from the simplicity of the model: considering a state space $N$, finite or infinite configurations are finite or countable sets of values in $N.$ This is why we begin with giving a short description of this setting, and describe a differentiable structure that can fit with easy problems of dynamics. This structure, which can be seen either as a Fr\"olicher structure or as a diffeological one, is carefully described and the links between these two frameworks are summarized in the appendix. We also give a result that seems forgotten in the past literature: the infinite configuration space used in e.g. \cite{Alb} is an infinite dimensional manifold. 

But the main goal of this paper is to include one dimensional turbulence effects (in particular period doubling, see e.g. \cite{Hil,Lesne}) in the dynamics described by finite or infinite configuration spaces. For this, we change the metric into the Hausdorff metric. This enables to ``glue together'' two configurations into another one and to describe shocks. 
Therefore, the dynamics on this modified configuration space are described by multivalued paths that are particular cases of graphs on $N,$ which explains the terminology: ``branched configuration''. 

We finish with the description of the configurations used e.g. in optimal transport, the space of probability measures, and show how they can also furnish configurations for uncompressible fluids. In these settings, branching effects are well-known and sometimes obvious, and we do not need any adaptation of the framework to obtain a full description of them.  Therefore, what we call branched configurations appears as an intermediate (an we hope useful) step between dynamics of e.g. a N-body problem and e.g. wave dynamics.    
\section{Configuration spaces} \label{gedebut}

We first recall the definition of configuration spaces
\subsection{Dirac configuration spaces on a locally compact manifold} \label{4.1}
Let us describe step by step a way to build infinite configurations, as they are built in the mathematical literature. 
We explain each step with the configurations already defined in e.g. \cite{Alb} and \cite{Str}, the generalization will be discussed later in this paper. 
A set of $1$-configurations is a set $\Gamma^1$ of objects that are modelizations of physical quantities. 
For example, in the settings \cite{Alb, Str, FKK0, FKK1, FKK2, FKK3, HKPR, Ism},  the physical quantity modeled is the position of one particle.
The whole world is modeled as a locally compact manifold $N$, and the set of 
$1$-configurations is itself $N$, or equivalently 
the set of all Dirac measures on $N$.

Let $I$ be a set of indexes. $I$ can be countable or uncountable.
 We define the \textbf{indexed} 
(or the \textbf{ordered} if $I \subset \N$ equipped with its total order)
 configuration spaces. For this, we need to define a 
symmetric binary relation relation $\U$
on $\Gamma^1$, that express the compatibility of two physical quantities.
We assume also that $\U$ has the following property:
\begin{equation}\label{incomp}
\forall (u,v) \in (\Gamma^1)^2, \quad u \U v \Rightarrow u \neq v.
\end{equation}
In the settings \cite{Ism}, \cite{Alb} 
and \cite{Str}, two particles 
cannot have the same position. Then, for $x, y \in N$, 
$$ x \U y \Longleftrightarrow x \neq y.$$ 
\vskip 10pt \noindent
With these restrictions, we can define the indexed or ordered configuration spaces :

\begin{eqnarray*} O\Gamma^n & = & \lbrace (u_1, ... , u_n) \in (\Gamma^1)^n 
\hbox{ such that, if } i \neq j, \quad u_i \U u_j \rbrace \\
 O\Gamma &=&  \coprod_{n \in \N^*} O\Gamma^n \hbox{ and}\\
O\Gamma^I &=& \lbrace (u_n)_{n \in I} \in (\Gamma^1)^I 
\hbox{ such that, if } i \neq j, \quad u_i \U u_j \rbrace.\end{eqnarray*}

The general configuration spaces are not ordered. Let $\Sigma_n$ be 
the group of bijections on $\N_n$, and $\Sigma_I$ be the set of bijections on $I$. We 
can define two actions: 
\begin{eqnarray*} \Sigma_n \times O\Gamma^n & \rightarrow & O\Gamma^n \\
                              (\sigma, (u_1, ... , u_n) ) &   \mapsto & (u_{\sigma(1)}, ... , u_{\sigma(n)}) \end{eqnarray*}

and its infinite analog: 
\begin{eqnarray*} \Sigma_I \times O\Gamma^I & \rightarrow & O\Gamma^I \\
                              (\sigma, (u_n)_{n \in I} ) &   \mapsto & (u_{\sigma(n)})_{n \in I} \end{eqnarray*}

\noindent
where $\Sigma_I$ is a subgroup of the group of bijections of $I.$ 
In the sequel, $I$ is  countable with discrete topology,which avoids 
topological problems on $\Sigma_I$ as in more complex examples. 
Then, we define general configuration spaces:

\begin{eqnarray*} \Gamma^n &=& O\Gamma^n  / \Sigma_n \\
 \Gamma &=&   \coprod_{n \in \N^*} \Gamma^n\\
\Gamma^I &= &O\Gamma^I  / \Sigma_I \end{eqnarray*}

\subsection{Configuration spaces on more general settings} \label{4.2}
In the machinery of the last section, the properties of the base manifold$N$ are not used in the definition of the space $\Gamma^1$. This is why the starting point can be $\Gamma^1$ instead of $N$, and we can give it the most general differentiable structure. 
Let us first consider the most general case: 
\begin{Proposition}
If $\Gamma^1$ is a diffeological space, then $\Gamma^n$, $\Gamma$ and $\Gamma^I$ are diffeological spaces
\end{Proposition}

\textbf{Proof.} $(\Gamma^1)^n$ (resp. $(\Gamma^1)^\N$) is a diffeological space
 according to Proposition \ref{prod1} (resp. Proposition \ref{prod3}), so that, 
$O\Gamma^n$ (resp. $O\Gamma^I$)  is a diffeological space as a subset of 
$(\Gamma^1)^n$ (resp. $(\Gamma^1)^\N$). Thus, $\Gamma^n = O\Gamma^n / \Sigma_n$ 
(resp. $\Gamma^I = O\Gamma^I / \Sigma)$ has the quotient 
diffeology by Proposition \ref{quotient}, which ends the proof. \qed

\vskip 12pt
Let us now turn to the cases where $\Gamma^1$ has a stronger structure. 
We already know that $\Gamma^n$ is a manifold if $\Gamma^1$ is a manifold. 

\begin{Proposition} \label{frolicherfini}
If  $\Gamma^1$ is a Fr\"olicher space, then $\Gamma^n$ (and hence $\Gamma$) is  a Fr\"olicher space.
\end{Proposition}

\textbf{Proof.}
Adapting the last proof, using Proposition \ref{prod2} instead of Proposition 
\ref{prod1}, we get that $O\Gamma^n$ is a Fr\"olicher space if $\Gamma^1$ is a 
Fr\"olicher space. Let us now build a generating set of functions for the 
Fr\"olicher structure on $O\Gamma^n$. Let $f: O\Gamma^n \rightarrow \R$ be a smooth 
map. We define the symmetrization of $f$: $$\tilde f : (u_1,... , u_n) \in O\Gamma^n 
\mapsto \tilde f(u_1,...,u_n) = {1 \over n!} \sum_{\sigma \in   \Sigma_n} f(u_{\sigma
(1)}, ... ,u_{\sigma(n)} ).$$
The set of functions $\{ \tilde f \}$ generate the contours on $O\Gamma^n$, and 
then, passing to the quotient, generate the contours on $\Gamma^n$. \qed
 
Setting a Fr\"olicher structure on indexed configurations is only a straightforward 
consequence of Proposition \ref{prod3}:

\begin{Proposition}
If  $\Gamma^1$ is a Fr\"olicher space, then $O\Gamma^I$ is  a Fr\"olicher space.
\end{Proposition}

But the problem of a Fr\"olicher structure on $\Gamma^I$ a little bit 
more complicated; 
let us explain why and give step by step the construction of the Fr\"olicher 
structure.
We first notice that the contours of the Fr\"olicher push forward naturally 
by the quotient map $O\Gamma^I \rightarrow \Gamma^I.$ But, if one wants to describe 
a generating set of functions, by Proposition \ref{prod2}, one has to
consider all combinations of a finite number of smooth functions 
$\Gamma^1 \rightarrow \R.$  This generating set does not contain any 
$\Sigma^I$-invariant function, except constant functions. This is why the approach 
used in the proof of the Proposition \ref{frolicherfini} cannot be applied here.
For this, one has to consider the set
$$\F_{eq} = \{ f : O\Gamma^I \rightarrow \R \in \F \hbox{ such that } 
f \hbox{ is } \Sigma^I \hbox{-invariant } \}$$
of equivariant functions on  $O\Gamma^I.$ This discusssion can become very quickly na\"ive and we prefer to leave this question to more applied works in order to fit with known examples instead of dealing with 
too abstract considerations.
\subsection{$\Gamma^{t,\infty}:$ an infinite configuration space with manifold structure}

If $\Gamma^1$ is endowed with a topology, which is 
authomaticaly the case when it is a 
diffeological space, $\Gamma^I$ is endowed with the trace of the product topology. 
But, for many reasons, one can want to restrict to some type of infinite 
configurations. Here, the diffeology (or the topology) endowed is still the trace diffeology (or topology), 
but one can modify the definition of $\Gamma^I$ in order to take under consideration 
the topology of $\Gamma^1$ in order to preserve physical realism. 
 
Fro example, if we want to recover the example of \cite{Alb}, 
which was one of our starting points, we have to remark that
 $\Gamma^I$ can be seen as the set of countable subsets of $\Gamma^1 = N$, 
which means that a subset $u \in \Gamma^I$ can have some accumulation point.  
In that case, the topology of $\Gamma^1$, and more precisely the sequential convergence in $\Gamma^1$, 
is taken into account for the \underline{definition} of another infinite configuration space, 
and not only for its topological structure. This leads to the following definitions, with $\Gamma^1 = N$:
$$O\Gamma^{t,\infty} = \lbrace (u_n)_{n \in \N} \in O\Gamma^I \hbox{ such that } \forall 
\hbox{ compact subset } K, \quad |\{u_n; n \in \N \}\cap K|<\infty \rbrace$$
and   
 $$\Gamma^{t,\infty} = \lbrace u\in \Gamma^I \hbox{ such that } \forall 
\hbox{ compact subset } K, \quad |u \cap K|<\infty \rbrace.$$
We now recognize the announced example. 
When $\Gamma^1$ is a locally compact manifold $N$, 
the corresponding infinite configuration space
$\Gamma^{t,\infty}$ is the one described in \cite{Alb}. 
In this case, since $\Gamma^1$ is a locally compact Lindel\"off manifold, 
we can build on $\Gamma^{t,\infty}$ a Fr\"olicher structure. 
After this review, we remark the following:

\begin{Lemma}
Let $u \in O\Gamma^{t,\infty}.$ We define 
$$ \varphi_u : u \rightarrow \R_+$$
by $$\varphi_u(x) = \inf\{d(x,x')| x' \in u \wedge x \neq x'\}.$$
Then, $\forall u \in O\Gamma^{t,\infty},$ $\varphi_u$ is $\R_+^*-$valued.   
\end{Lemma}

\noindent
\textbf{Proof.}
Let $U$ be a relatively compact neighborhood of $x.$ Then $\bar{U} \cap u$
is a finite set, so that
$$ \varphi_u(x) = \inf\{d(x,x')| x' \in u \wedge x \neq x'\} = \inf\{d(x,x')| x' \in \bar{U} \cap u \wedge x \neq x'\} > 0.$$ \qed

Let $p = \dim N$ and let $B_\epsilon$ be the Euclidian ball of $\R^p$ of radius $\epsilon.$ Let us now fix $\{\phi_\alpha : B_1 \rightarrow N\}_{\alpha \in \Lambda}$ a smooth atlas on $N.$ Let $u \in O\Gamma^{t,\infty}$ and let $x \in u.$ We fix $\alpha_x\in \Lambda$ an index such that $x \in \phi_{\alpha_x}(B_1).$ By translation in $\R^p,$ we change $\phi_{\alpha_x}$ into a map $\tilde\phi_{\alpha_x}$ such that $\tilde\phi_{\alpha_x}(0)=x.$ Let 
$$\epsilon_x = \frac{1}{2} \sup \{ \epsilon >0 | d(\tilde\phi_{\alpha_x}(B_\epsilon), u-\{x\})>0\}$$
This non-zero number exists by the previous lemma. 
Then we define
$ \Phi_{x,u} : B_1 \rightarrow N$
by $$\Phi_{x,u}(z) = \tilde\phi_{\alpha_x} (\epsilon_x z).$$
By construction,  $\Phi_{x,u}(B_1) \cap u = \{x\}.$
Now we can state the theorem:
\begin{Theorem}
Let $u \in O\Gamma^{t,\infty}.$ We define 
$$ \Phi_u : B_1^\N \rightarrow N^\N$$
by $$\Phi_u(z_0,z_1,...) = \left( \Phi_{u_0,u}(z_0),  \Phi_{u_1,u}(z_1),...\right).$$
Then, 
\begin{itemize}
\item $\Phi_u ( B_1^\N) \subset O\Gamma^{t,\infty};$
\item Let $(u,u') \in (O\Gamma^{t,\infty})^2$ If $\Phi_u(B_1^\N) \cap \Phi_{u'}(B_1^\N) \neq \emptyset,$ the map $\Phi_u^{-1} \circ \Phi_{u'}$ is smooth (in the Fr\"olicher sense) from $\Phi_{u'}^{-1} \left( \Phi_{u} (B_1^\N) \right)$ to $\Phi_{u}^{-1} \left( \Phi_{u'} (B_1^\N) \right).$
\end{itemize}

\end{Theorem}

\noindent
\textbf{Proof.} The proof is obvious, by construction of the maps $\Phi_u.$ \qed

Thus, we get:
\begin{Cor} if we equip $B_1^\N$ with a metric induced by a norm $||.||$  such that $||.|| < ||.||_{l^\infty},$ we give to $O\Gamma^{t,\infty}$ a structure of Banach manifold with atlas $\left\{ \Phi_u | u \in O\Gamma^{t,\infty} \right\}.$
\end{Cor}
 And by projection $O\Gamma^{t,\infty} \rightarrow \Gamma^{t,\infty},$ we also get an atlas on (non idexed) configurations $\Gamma^{t,\infty}.$

\section{``Topological'' configuration spaces}

In this section, we present examples of 1-configurations and 
their associated configuration spaces. 
Manifolds will replace the Dirac measures used in \cite{Alb}. 
In the sequel, $N$ is a Riemannian smooth locally compact manifold.
the 1-configurations considered keep their \textbf{topological} properties, as in 
the model of elastodymanics (see e.g. \cite{Marsden}) or in various quantum field theories.
Notice also that we do not give compatibility conditions between two 
1-configurations: we would like to give the more appropriate conditions in order to fit with the applied models,  
this is why we leave this point to more specialized works. 
   
\subsection{Topological 1-configurations}

We follow here, for example, \cite{Marsden}. 

\begin{Definition}
Let $M$ be a smooth compact manifold and $N$ an arbitrary manifold. 
We set 
$$\Gamma^1_{e}(M,N) = C^\infty(M,N).$$
\end{Definition}

One can also only consider embeddings, and set : 
$$\Gamma^1_{m}(M,N) = Emb(M,N),$$
where $dim M < dim N$, and $Emb$ is the set of embeddings. The things run as in the 
first case, since $Emb(M,N) \subset C^\infty(M,N)$ is an open subset 
of $C^\infty(M,N).$

\subsection{Examples of topological configurations}
\subsubsection{Links}
Let $\Gamma^1 = Emb(S^1,N)$
Here, we fix the uncompatibility relation as 
$$ \gamma \U \gamma' \Leftrightarrow \gamma(S^1) \cap \gamma'(S^1)\neq \emptyset.$$
Then, $\Gamma^n_{link}$ is the space of $n-$links of class $C^k$, whick is a Fr\'echet manifold.
\subsubsection{Triangulations} 
Consider the n-simplex 
$$\Delta_n = \left\{ (t_0, ..., t_n) \in \R^{n+1} | \sum_{i=0}^n t_i = 1 \right\}.$$
If $N$ is a n-dimensional manifold, a (finite) triangulation $\sigma$ of $N$ is such that:
\begin{enumerate}
\item $\sigma \in \left(\Gamma^1_{m}(\Delta_n,N)\right)^{|\sigma|}$
\item \label{compatible} Let $(\tau, \tau') \in \sigma^2 $ such that $\tau \neq \tau',$ then $Im(\tau) \cap Im(\tau') $ is a simplex or a collection of simplexes of each border
$Im(\partial \tau)$ and $Im(\partial\tau').$
\item $$\bigcup_{\tau \in \sigma} Im \tau = N.$$
\end{enumerate}

We get by condition \ref{compatible} a compatibility condition $\U,$ for which we can build  $O\Gamma(\Delta_n, N)$, $\Gamma(\Delta_n,N), $ $O\Gamma^\infty(\Delta_n, N)$ and $\Gamma^\infty(\Delta_n,N). $ If $N$ is compact, the set of triagulations of $N$ is a subset of $\Gamma(\Delta_p,N).$ If $N$ is non compact and locally compact, the set of triangulations of $N$  is a subset of $\Gamma^\infty(\Delta_n,N).$

More generally, for $ p \leq n,$ one can build $O\Gamma(\Delta_p, N)$, $\Gamma(\Delta_p,N), $ $O\Gamma^\infty(\Delta_p, N)$ and $\Gamma^\infty(\Delta_p,N). $ This example will be discussed in the section \ref{br}.

\subsubsection{Strings and membranes}
A string is a smooth surface $\Sigma,$ possibly with boundary, embedded in $\R^{26}.$  A membrane is a manifold $M$ of higher dimension embedded in some $\R^k.$  We recover here some spaces of the type $\Gamma^1_m$, which will be also discussed in section \ref{br}.

\section{Branched configuration spaces} \label{br}

\subsection{Dirac branched configurations} \label{brnchedcf}
As we can see in section \ref{4.1}, finite configurations $\Gamma$ are made of a countable disjoint union.
We now fix a metric $d$ on $N.$
The idea of branched configurations is to glue together the components 
$\Gamma^n$ on the generalized diagonal, namely, we define the following distance 
on $\Gamma$:   
\begin{Definition}
Let $(u,v)\in \Gamma^2.$ $$d_\Gamma(u,v) = \sup_{(x,x') \in u\times v} \{ d(x,v), d(x',u)\} $$
\end{Definition}

\begin{Proposition}
$d_\Gamma$ is a metric on $\Gamma .$
\end{Proposition}
\noindent
\textbf{Proof.}
We remark that $d_\Gamma$ is the Hausdorff distance restricted to $\Gamma.$ \qed

%Now, we remark that $f \in C^\infty(M,\R)$ induces a map on $\Gamma$ by $$f(u) = \frac{1}{|u|} \sum_{x \in u} f(x).$$ 
% Let $\F_0 = C^\infty(M,\R).$ Following \cite{KM}, the set of paths 
%$\C = \{c \in C^0(\R, \Gamma) | 
%\forall f \in \F_0, f \circ c \in C^\infty(\R,\R)\}$
%defines a Fr\"olicher space (\cite{FK}, see e.g. \cite{KM}). 
The following proposition traduces the change of topology of $\Gamma$ into $\B\Gamma$ by cut-and paste property:

\begin{Proposition}
$\forall n \in \N^*, \B\Gamma^{n+1} = \Gamma^{n+1} \coprod_{\B\Gamma^n} \B\Gamma^n$
 where the identification is made along the trace on $\Gamma^{n+1}$ of the $d_\Gamma-$neighborhoods of $\B\Gamma^n$ in $N^{n+1} \supset O\Gamma^{n+1}.$
\end{Proposition}

We remark that we can also define a Fr\"olicher structure on $\B\Gamma$ with generating set of functions the set 
$$\left\{ u \in \O\Gamma \mapsto \frac{1}{|u|}\sum_{x\in u} f(x) | f \in C^\infty(N,\R \right\}.$$
This structure will be recovered later in this paper.

\subsection{Examples of topological branched configurations}

\subsubsection{The path space, branched paths and graphs}

Let $$\Gamma^1_e([0;1];N) = C^\infty([0;1];N)$$ be the space of smooth paths on $N.$
A path $\gamma$ has a natural orientation, and has a beginning $\alpha(\gamma)$ and an end $\omega(\gamma).$ We define a compatibility condition 
$$ \gamma \U \gamma' \Leftrightarrow \left( (Im\gamma \neq Im \gamma') \vee \left(\alpha(\gamma), \beta(\gamma)) \neq (\alpha(\gamma'), \beta(\gamma')\right)  \right)$$
and we remark that the set of piecewise smooth paths on $N$ is a subset of $O\Gamma_e([0;1];N),$ saying that $(\gamma_1,... \gamma_p) \in O\Gamma_e^p([0;1];N)$ is a piecewise smooth path if and only if $$\forall i \in \N_{p-1}, \omega(\gamma_i) = \alpha(\gamma_{i+1}).$$ This relation, stated from the natural definition of the composition $*$ of the groupo\"id of paths, is not unique and can be generalized.
\begin{Definition}
 Let $(\gamma_1,\gamma_2) \in O\Gamma^2_e([0;1];N)$
and let $\gamma_3 \in C^\infty([0;1];N).$ We define the equivalence relation $\sim$ by $$(\gamma_1,\gamma_2) \sim \gamma_3 \Leftrightarrow \gamma_3 = \gamma_2 * \gamma_1.$$
 \end{Definition}

The maps $\alpha: \gamma \mapsto \alpha(\gamma)$ and $\omega: \gamma \mapsto \omega(\gamma)$ extends to ``set theorical'' maps $O\Gamma_e^k([0;1];N) \rightarrow N^k$ and $\Gamma_e([0;1];N) \rightarrow \Gamma(N).$
The following is now natural:

\begin{Definition}
A \textbf{branched path} is an element $\gamma$ of $O\Gamma\left(O\Gamma_e([0;1];N) / \sim \right)$ such that, if $\gamma \in O\Gamma^k\left(O\Gamma_e([0;1];N) / \sim \right),$ $$\forall i \in \N_{k-1},  \omega(\gamma_i)  = \alpha(\gamma_{i+1}) \hbox{ in } \Gamma(N).$$
\end{Definition}

\noindent
\textbf{Example.}
Let us consider the following paths $[0;1] \rightarrow \R^2:$ 
\begin{itemize}
\item $\gamma_1(t) = (t-2;0)$
\item $\gamma_2(t) = (cos(\pi (1-t)); sin(\pi t))$
\item $\gamma_3(t) = (cos(\pi (1-t)); -sin(\pi t))$
\item $\gamma_4(t) = (t+1;0)$
\end{itemize}
Then, $$\omega(\gamma_1) = \alpha(\gamma_2) = \alpha(\gamma_3)$$
and $$\omega(\gamma_2) = \omega(\gamma_3) = \alpha(\gamma_4).$$
This shows that $$\left( \gamma_1, (\gamma_2, \gamma_3) , \gamma_4\right) \in O\Gamma^3(O\Gamma_e([0;1];\R^2)$$ is a branched path of $\R^2.$

\subsubsection{Alternate approach to branched paths: branched sections of a fiber bundle}
Let $\pi:F\rightarrow M$ be a fiber bundle of typical fiber $F_0.$

Here, $n \in \N^* \cup \infty .$
Let $\pi:F \rightarrow M$ be a fiber bundle over $M$ with typical fiber $F_0.$
Let $$\Gamma_M^n(F) = \left\{ u \in \Gamma^n(F)| |\pi(u)| = 1\right\}.$$
This is trivially a fiber bundle of basis $M$ with typical fiber $\Gamma^n(F_0).$ 
\begin{Definition}
A \textbf{non-section} of $F$ is a section of $\Gamma_M^n(F)$ which cannot be decomposed into $n$ sections of $F.$
\end{Definition}
We define also $\Gamma_M(F) = \coprod_{n \in \N^*}\Gamma^n_M(F), $ and also $\Gamma^I_M(F)$ the non sections based of $\Gamma^I(F_0).$
We can  define the same way $\B  \Gamma_M(F)$ using the branxhed configuration space intead of the configuration space, since the definitions from the  set-theoric viewpoint are the same. 
\vskip 12pt
\noindent
\textbf{(Toy) Example.} Let us consider the following example: $X = \R^3 \times \{ \hbox{up ; down} \}$, and 
$\Gamma^1(X)=\R^3 \times \{ \{ \hbox{up} \} ;
\{ \hbox{down} \} ; \{ \hbox{up ; down} \}\},$ that models 
the position of an electron in the 3-dimensional space $\R^3$, associated to its spin. 
When the electron spin cannot be determined (i.e. out of the action of adequate electromagnetic fields), 
the picture proposed by Schr\"odinger is to consider that its spin is both up and down 
(this picture is also called 
the ``Schr\"odinger cat'' when we replace ``up'' and ``down'' by ``dead'' and ``alive'').  

\vskip 12pt
Let us now consider the Fr\"olicher structure described on section \ref{brnchedcf}. It is based on the natural diffeology carried by each $\Gamma^n(F_0)$ ($n \in \N^*$) and by the set of paths $\p'_1$  that are paths $\gamma : \R \rightarrow \Gamma(F_0)$ such that $\exists (m,n) \in (\N^*)^2,$

\begin{itemize}
\item $\gamma|_{]-\infty;0]} $ is a smooth path on $\Gamma^m(X),$ 
\item  $\gamma|_{]0; +\infty;0[} $ is a smooth path on $\Gamma^n(X),$
\item Let $l \in \gamma(0).$ Then for any smooth map $ f: F \rightarrow \R,$ the infinite jet of $$\sum f \circ c_-$$ where $c_-$ are  
the trajectories going to l in $0^-$ equals to the sum of the infinite jet of $$\sum f \circ c_+$$ where $c_+$ are
 the trajectories coming from $l$ in $0^+.$  

\end{itemize}
Remark that the last condition comes from the smoothness required for each map $f \circ \gamma,$ with $f \in C^\infty(F_0,\R).$ This fits with the (fiberwise) fr\"olicher structure of $\B\Gamma(\F_0).$ Then, a finitely branched section of $F$ is a smooth section of $\B\Gamma_M(F).$
The first examples that come to our mind are the well-known branched processes, and we can wonder some deterministic analogues replacing stochastic processes by dynamical systems. Let us here sketch a toy example extracted from the theory of turbulence:
\vskip 12pt
\noindent
\textbf{Example: equilibum of mayflies population}
Assuming that Mayflies live and die in the same portion of river, 
the population $p_{n+1}$ at the year $(n+1)$ is obtained from the population $p_n$ at the year $n$ (after normalization procedure) 
by the formula $$ p_{n+1} = A p_n(1-p_n)$$
where $A \in [0;4]$ is a constant coming from the environmental data. 
For $A$ small enough, the fixed point of the so-called ``logistic map'' $\phi_A(x)=Ax(1-x)$ is stable, hence the population $p_n$ tends to stabilize around this value. But when $A$ is increasing, the fixed point becomes unstable and $p_n$ tends to stabilize around $2^k$ multiple values which are the stable fixed points of of the map $\phi^{2^k}$ obtained by composition rule. 

Now, assume that we consider a river (or a lake), modelized by an interval 
(or an open subset of $\R^2$) that we denote by  $U$, where the parameter $A$ is a smooth map 
$U \rightarrow [0;4].$ The parameter $A$ is a smooth map $U \rightarrow [0;4]$ and the cardinality 
of configuration of equilibrum depends on the value of $A.$ 
\subsubsection{Non sections in higher dimensions}
The example of a lake where mayflies live and die gives us a nice example of branched surface viewed as an element of $\B\Gamma_{\R^2}(\R^2 \times [0;1]).$ The same procedure can be implemented in gluing simplexes, or strings or membranes along their borders to get branched objects, but we prefer to postpone this problem to a work in progress where links with stochastic objects should be performed.

\section{Measure-like configurations: an example at the borderline of branched configurations and dynamics on probability spaces}

Dymanics on probability space is a fast-growing subject and is shown to give ise to branched geodesics \cite{Ohta}. Following the same procedure as for branched topological configurations,we show her how a restricted space fits with particular goals. The goals described here are linked with image recognition for the configuration space $\Gamma_h^1$ above, and to umcompressible fluid dynamics when we equip   $\Gamma_h^1$ with the diffeology $\p_0$ of constant volume above.
Let $C^0_c(N)$ be the set of compactly supported $\R$-valued smooth maps on $N$. 
We define the relation of equivalence $\rel$ by:
$$ f \rel g \quad \Leftrightarrow \quad Supp(f) = Supp(g) .$$
Let us first give the definition of the set of 1-configurations :

\begin{Definition}
We set 
$$\Gamma^1_h(N)= C^0_c(N) / \rel .$$
\end{Definition}

Such a space is not a manifold, but we show that it carries a natural diffeology. 
$C^0_c(N)$ is a topological vector space, and hence carries a natural 
diffeology $\p_0$.
We define the following:

\begin{Definition} \label{vol}
Let $\p_1 \subset \p_0$ be the set of 
$\p_0$-plots $p : O \rightarrow C^0_c(N)$ such that,
for any open subset $A$ with compact closure $\overline{A}$ of $N$, 
for any open subset $O'$
 of $O$ such that $\overline{O'} \subset O$, 
if $$ p( O')(\overline{A} \backslash A) = {0} ,$$
the map 
$$ x \in O \mapsto Vol(Supp(p(x)) \cap A )$$
is constant on $O'$, where $Vol$ is the Riemannian volume. 
\end{Definition}
 
This technical condition ensures that the volume of any connected component 
of the support of 
$p(x)$ is constant.
$\p_1$ is obviously a diffeology on $C^0_c(N)$, and we can state

\begin{Proposition}
Let $\p =$ $ \p _1 / \rel$.
$(\Gamma^1_h(N), \p)$ is a diffeological space.
\end{Proposition}

The proof is a straightforward application of Definition \ref{quotient}. 
It seems difficult to give this space a structure of Fr\"olicher space, 
or even a natural topology except the topology of vague convergence of measures, 
which is not the topology induced by the diffeology we have defined. 
As a consequence, we can only state that 
the well-defined configuration spaces $\Gamma$ and $\Gamma^I$ are 
diffeological spaces. The technical conditions of Definition \ref{vol}
 ensures that volume preserving is a consequence of smoothness with respect to $\p_1,$ and hence is particularily designed for (vicou) uncompressible fluid dynamics. 
A $1-$conffuguration can have many connected components, and therefore branching effects are included in the definition of $1-$configurations. 
One could understand $n-$configurations as the presence of $n$ (non mixing) fluids.

 \section*{Appendix: Preliminaries on differentiable structures} \label{1.}
The objects of the category of -finite or infinite- dimensional smooth manifolds is
made of topological spaces $M$ equipped with a collection of charts called maximal atlas that
enables one to make differentiable calculus. But there are some examples
where a differential calculus is needed where as no atlas can be defined. To circumvent this problem,
several authors have independently developped some ways to define differentiation without defining charts.
We use here three of them. The first one is due to Souriau \cite{Sou},
the second one is due to  Sikorski, and the third one is a setting closer to the setting of differentiable manifolds
is due to Fr\"olicher (see e.g. \cite{CN} for an introduction on these two 
last notions). In this section, we review some basics on these three notions.

\subsection*{Souriau's diffeological spaces, Sikorski's differentiable spaces, Fr\"olicher spaces}

\begin{Definition}
Let $X$ be a set.

\noindent
$\bullet$ A \textbf{plot}	 of dimension $p$ (or $p$-plot) on $X$ is a map
from an open subset $O$ of $\R^p$ to $X$.

\noindent
$\bullet$ A \textbf{diffeology} on $X$ is a set $\p$ of plots on $X$ such that,							  
for all $p \in \N$,

- any constant map $\R^p \rightarrow X$ is in $\p$;

-  Let $I$ be an arbitrary set; let 
$\{f_i : O_i \rightarrow X\}_{i \in I} $ be a family of maps that extend to a map 
$f: \bigcup_{i \in I}O_i \rightarrow X$. If $\{f_i : O_i \rightarrow X\}_{i \in I} \subset \p $, 
then  $f \in \p$.

- (chain rule) Let $f\in \p$, defined on $O \subset \R^p$. Let $q \in \N$,
$O'$ an open subset of $\R^q$ and $g$ a smooth map (in the usual sense) from $O'$ to $O$.
Then, $ f \circ g \in \p$.

\vskip 6pt
\noindent
$\bullet$ If $\p$ is a diffeology $X$, $(X,\p)$
is called \textbf{diffeological space}.

\noindent Let $(X,\p)$ et $(X',\p')$ be two diffeological spaces,
a map $f : X \rightarrow X'$ is \textbf{differentiable} (=smooth)
if and only if $f \circ \p \subset \p'$.
\end{Definition}

\noindent
\textbf{Remark.} Notice that any diffeological space $(X,\p)$
can be endowed with the weaker topology such that all the maps that are in $\p$ are continuous. 
But we prefer to mention this only for memory as well as other questions 
that are not closely related to our construction, and stay closer to the goals of this paper.

Let us now define the Sikorski's differential spaces.
Let $X$ be a Haussdorf topological space.

\begin{Definition}

\noindent
$\bullet$ A (Sikorski's) \textbf{differential space} is a pair $(X, \F)$ where
$\F$ is a family of functions $ X \rightarrow \R$ such that

- the topology of $X$ is the initial topology with respect to $\F$

- for any $n \in \N$, for any smooth function $\varphi : \R^n \rightarrow \R$, for any
$(f_1,...,f_n) \in \F^n$, $\varphi \circ (f_1,...,f_n) \in \F$.

\noindent
$\bullet$ Let $(X,\F)$ et $(X',\F')$ be two differential spaces,
a map $f : X \rightarrow X'$ is \textbf{differentiable} (=smooth)
if and only if $ \F' \circ f \subset \F.$
\end{Definition}

We now introduce Fr\"olicher spaces.

\begin{Definition}
$\bullet$ A \textbf{Fr\"olicher} space is a triple $(X,\F,\C)$ such that

- $\C$ is a set of paths $\R \rightarrow X$,

- A function $f: X \rightarrow \R$ is in $\F$ if and only if for any $c \in \C$, $f\circ c \in C^\infty(\R,\R)$;

- A path $c: \R \rightarrow X$ is in $\C$ (i.e. is a \textbf{contour}) if and only if for any $f \in \F$, $f\circ c \in C^\infty(\R,\R)$.

\vskip 5pt
$\bullet$ Let $(X,\F,\C)$ et $(X',\F',\C')$ be two Fr\"olicher spaces,
a map $f : X \rightarrow X'$ is \textbf{differentiable} (=smooth)
if and only if $\F'\circ f \circ \C \in C^\infty(\R,\R)$.
\end{Definition}

Any family of maps $\F_g$ from $X$ to $\R$ generate a Fr\"olicher structure  $(X,\F,\C)$, 
setting \cite{KM}:

- $\C = \{ c : \R \rightarrow X \hbox{ such that } \F_g \circ c \subset C^\infty(\R,\R) \}$

- $\F = \{  f: X \rightarrow \R \hbox{ such that } f \circ \C \subset C^\infty(\R,\R) \}.$

One easily see that $\F_g \subset \F$. This notion will be useful in the sequel to describe
in a simple way a Fr\"olicher structure. 

A Fr\"olicher space, as a differential space, carries a natural topology, which is the pull-back topology of $\R$
via $\F$. In the case of a finite dimensional differentiable manifold, 
the underlying topology of the Fr\"olicher structure is the 
same as the manifold topology. In the infinite dimensional case, these two topologies differ very often.

\vskip 12pt
 
In the three previous settings, we call $X$ a \textbf{differentiable space}, 
omitting the structure considered. Notice that, in the three previous settings, 
the sets of differentiable maps between two differentiable spaces
satisfy the chain rule. 
Let us now compare these three settings:
One can see (see e.g. \cite{CN}) that we have the following, given at each step by forgetful functors:
\vskip 10pt
\begin{tabular}{ccccc}

smooth manifold & $\Rightarrow$ & Fr\"olicher space & $\Rightarrow$ & Sikorski differential space\\
\end{tabular}
\vskip 10pt
Moreover, one remarks easily from the definitions that, if $f$ is a map from a 
Fr\"olicher space $X$ to a Fr\"olicher space $X'$, $f$ is smooth in the sense of 
Fr\"olicher if and only if it is smooth in the sense of Sikorski. 

\vskip 12pt

One can remark, if $X$ is a Fr\"olicher space, 
we define a natural diffeology on $X$ by \cite{Ma}: 
$$\p(\F) = \coprod_{p \in \N} \{ \, f \hbox{ $p$-paramatrization on } X; \,
 \F \circ f \in C^\infty(O,\R) \quad \hbox{(in the usual sense)}\}.$$
With this construction, we get also a natural diffeology when $X$ 
is a Fr\"olicher space. In this case, one can easily show the following:

\begin{Proposition} \label{Frodiff} \cite{Ma}
Let $(X,\F,\C)$ and $(X',\F ', \C ')$ be two Fr\"olicher spaces. 
A map $f : X \rightarrow X'$ is smooth in the sense of  Fr\"olicher 
if and only if it is smooth for the underlying diffeologies.
\end{Proposition}

Thus, we can also state: 
 
 \begin{tabular}{ccccc}

smooth manifold & $\Rightarrow$ & Fr\"olicher space & $\Rightarrow$ & Diffeological space\\
\end{tabular}

\subsection*{Cartesian products}

The category of Sikorski differential spaces is not cartesianly closed, 
see e.g. \cite{CN}. This is why we prefer to avoid the questions
related to cartesian products on differential spaces in this text, and focus 
on Fr\"olicher and diffeological spaces, since the cartesian product 
is a tool essential for the definition of configuration spaces. 

In the case of diffeological spaces, 
we have the following \cite{Sou}:

\begin{Proposition} \label{prod1} Let $(X,\p)$ and $(X',\p')$ be two diffeological spaces. We call
\textbf{product diffeology} on $X \times X'$ the diffeology $\p \times \p'$ made of plots
$g : O \rightarrow X \times X'$ that decompose as $g = f \times f'$,
where $f : O \rightarrow X \in \p$ and $f': O \rightarrow X' \in \p'$.
\end{Proposition}

Then, in the case of a Fr\"olicher space, we derive very easily, compare with e.g. \cite{KM}:
\begin{Proposition} \label{prod2}
Let $(X,\F,\C)$ and $(X',\F ', \C ')$ be two Fr\"olicher spaces, 
with natural diffeologies $\p$ and $\p '$ . There is a natural structure of Fr\"olicher space on 
$X \times X'$ which contours $\C \times \C'$ are the 1-plots of $\p \times \p'$.
\end{Proposition}

We can even state the following results in the case of infinite products. 

\begin{Proposition} \label{prod3} 
Let $I$ be an infinite set of indexes, that can be uncoutable.

$\bullet$ (adapted from \cite{Sou} )
Let $\lbrace (X_i,\p_i) \rbrace_{i \in I}$ be a family of diffeological spaces indexed by $I$. We call
\textbf{product diffeology} on $\prod_{i \in I} X_i $ the diffeology $\prod_{i \in I}\p_i $ made of plots
$g : O \rightarrow \prod_{i \in I}X_i$ that decompose as $g = \prod_{i \in I} f_i$,
where $f_i \in \p_i$. This is the biggest diffeology for which the natural projections are smooth.

$\bullet$ Let $\lbrace (X_i,\F_i,\C_i) \rbrace_{i \in I}$ be a family of Fr\"olicher spaces indexed by $I$, 
with natural diffeologies $\p_i$. There is a natural structure of Fr\"olicher space  
$( \prod_{i \in I}X_i, \prod_{i \in I}\F_i, \prod_{i \in I}\C_i)$ 
which contours $\prod_{i \in I} \C_i$ are the 1-plots of $\prod_{i \in I}\p_i$.
A generating set of functions for this Fr\"olicher space is the set of maps of the type:
$$ \varphi \circ \prod_{j \in J} f_j$$
where $J$ is a finite subset of $I$ and $\varphi$ is a linear map $\R^{|J|} \rightarrow \R$.
\end{Proposition}

\noindent
\textbf{Proof.}

$\bullet$ By definition, following \cite{Sou}, $\prod_{i \in I}\p_i $ is the biggest diffeology 
for which natural projections are smooth. Let $g: O \rightarrow X_i$ be a plot. 
$$g \in \prod_{i \in I}\p_i\Longleftrightarrow p_i \circ g \in \p_i, $$
where $p_i$ is the natural projection onto $X_i$, which gets the result.

$\bullet$ With the previous point and Proposition \ref{Frodiff}, 
we get the family of contours of the product Fr\"olicher space. \qed

\subsection*{Push-forward, quotient and trace}

We give here only the results that will be used in the sequel.

\begin{Proposition} \cite{Sou}
 Let $(X,\p)$ be a diffeological space, and let $X'$ be a set. Let $f:X \rightarrow X'$ be a surjective map. Then, the set
$$ f(\p) =\{ u  \hbox{ such that } u \hbox{restricts to some maps of the type } f \circ p;  p \in \p \} $$
is a diffeology on $X'$, called the \textbf{push-forward diffeology} on $X'$ by $f$.
\end{Proposition}

We have now the tools needed to describe the diffeology on a quotient:

\begin{Proposition} \label{quotient}
let $(X,\p)$ b a diffeological space and $\rel$ an equivalence relation on $X$. Then, there is 
a natural diffeology on $X / \rel$, noted by $\p / \rel$, 
defined as the push-forward diffeology on $X / \rel$ by the quotient 
projection $X \rightarrow X / \rel$.   
\end{Proposition}

Given a subset $X_0 \subset X$, where $X$ is a Fr\"olicher space or a diffeological space, we can define on trace structure on $X_0$, induced by $X$.

$\bullet$ If $X$ is equipped with a  diffeology $\p$, we can define a diffeology $\p_0$ on $X_0$ setting 
$$ \p_0 = \lbrace p \in \p \hbox{ such that the image of } p \hbox{ is a subset of  } X_0 \rbrace.$$

$\bullet$ If $(X, \F, \C)$ is a Fr\"olicher space, we take as a generating set of maps $\F_g$ on $X_0$
the restrictions of the maps $f \in \F$. In that case, the contours (resp. the induced diffeology) on $X_0$
are the contours (resp. the plots) on $X$ which image is a subset of $X_0$.

\end{document}